\documentclass[10pt]{amsart}
\usepackage[matrix,arrow,curve]{xy}
\usepackage{amssymb, amsthm}
\usepackage{supertabular}
\usepackage{amsfonts}
\usepackage{latexsym}
\usepackage{enumerate}
\usepackage[dvips]{graphics}

\usepackage{diagrams}
\usepackage{mathrsfs}
\author{Michael Farber and Mark Grant}
\address{Department of Mathematical Sciences, Durham University, South
  Road, Durham DH1 3LE}
\email{michael.farber@durham.ac.uk}
\email{mark.grant@durham.ac.uk}
\date{}

\title[Weights of cohomology classes]{Robot motion planning, weights of cohomology classes, and cohomology operations}

\newcommand{\bdm}{\begin{displaymath}}
\newcommand{\edm}{\end{displaymath}}

\newcommand{\Z}{\mathbf{Z}}
\newcommand{\co}{\colon\thinspace}

\newcommand{\TC}{\mathrm {TC}}

\newcommand{\wgt}{\mathrm{wgt}}
\newcommand{\swgt}{\mathrm{swgt}}

\newcommand{\genus}{\mathfrak{genus}}
\newcommand{\cat}{\mathrm{cat}}
\newcommand{\kk}{\mathbf k}
\newtheorem{thm}{Theorem}

\newtheorem{prop}[thm]{Proposition}
\newtheorem{definition}[thm]{Definition}

\newtheorem{lemma}[thm]{Lemma}

\newtheorem{cor}[thm]{Corollary}
\newtheorem{remark}[thm]{Remark}
\thanks{The authors were supported by a grant from the UK Engineering and Physical Sciences Research Council; the first author was also supported by a grant from the Royal Society.}


\begin{document}
\begin{abstract} The complexity of algorithms solving the motion planning problem is measured by a homotopy invariant $\TC(X)$ of the configuration space $X$ of the system. Previously known lower bounds for $\TC(X)$ use the structure of the cohomology algebra of $X$.
In this paper we show how cohomology operations can be used
  to shar\-pen these lower bounds for $\TC(X)$.
As an application of this technique we calculate explicitly the topological complexity
  of various lens spaces. The results of the paper were inspired by the work of E. Fadell and S. Husseini on
weights of cohomology classes appearing in the classical lower bounds for the Lusternik - Schnirelmann category.
In the appendix to this paper we give a very short proof of a generalized version of their result.
\end{abstract}

\maketitle

\section{Introduction}

The motion planning problem is a central theme of robotics. Given a mechanical system $S$,
a motion planning algorithm for $S$ is a function which associates to any pair of states $(A, B)$ of $S$ a continuous motion of the system starting at $A$ and ending at $B$.
If $X$ denotes the configuration space of the system, one considers the path fibration
\begin{eqnarray}\label{path}\pi\co X^I\to X\times X,\quad \pi(\gamma)=(\gamma(0),\gamma(1)), \quad \gamma: I\to X,\end{eqnarray}
and, in these terms, a motion planning algorithm for $S$ is a section (not necessarily continuous) of $\pi$.
The topological
  complexity of $X$, denoted $\TC(X)$, is defined to be the
genus, in the sense of Schwarz, of fibration (\ref{path}).
The concept $\TC(X)$ was introduced and studied in \cite{Far03}, \cite{Far04}; it is a measure of complexity of the problem of finding a
motion planning algorithm for a system whose configuration space is homotopy equivalent to $X$. A recent survey of related results can be found in \cite{Far06}.

A lower bound for $\TC(X)$ was given in \cite{Far03} in terms
of zero-divisors in cohomology with coefficients in a field $\kk$.
The
cup product map \begin{eqnarray}\label{iso}
\cup\co H^*(X;\kk)\otimes
H^*(X;\kk)\to H^*(X;\kk)\end{eqnarray}
is an algebra homomorphism, whose kernel is called the {\em
  ideal of zero-divisors}.
  The multiplicative structure on the left in (\ref{iso})
   is given by the formula $(\alpha\otimes \beta)(\gamma\otimes
\delta)=(-1)^{|\beta||\gamma|}\alpha\gamma\otimes\beta\delta.$
 A theorem from \cite{Far03} claims that {\it $\TC(X)$ is greater than the zero-divisors cup-length of $X$,} where
the latter is defined as the length of
the longest non-trivial product of elements in the
ideal of zero-divisors. This lower bound is sharp in many cases, and is easy to apply
as it only requires knowledge of the cohomology algebra of $X$.

In this paper we show how the lower bound for $\TC(X)$ mentioned above may
be improved upon using cohomology operations. We employ the notion of {\em weight} of a cohomology class with respect to a fibration,
introduced in \cite{FG06}, which generalises the notion of {\em category weight} developed by Fadell
and Husseini \cite{FH92} for estimating the Lusternik-Schnirelmann
category of a space.

The main result of this paper is described in \S \ref{sectcweight}. Here we introduce the notion of {\it excess} of a stable cohomology operation.
With the aid of this notion we find indecomposable zero-divisors having weight at least 2 with respect to the path fibration (\ref{path}).
This result is applied in \S \ref{seclens} to the computation of
the topological complexity of some lens spaces. In section \S \ref{secfibre} we give an upper bound for the topological complexity of a fibration which is
used in \S \ref{seclens}.

Fadell and Husseini showed that certain
cohomology classes
which are the images of Steenrod operations have weight at
least 2. As Appendix A to this paper we give a short proof of a more general result. Our Theorem \ref{FH} is more general than the result of \cite{FH92} in two respects: we allow more general cohomology operations, and give estimates for the strict category weight of Rudyak \cite{Rud99} as opposed to the usual category weight of Fadell and Husseini \cite{FH92}.

{\em All topological spaces are assumed to be path-connected, and all maps are continuous. All rings are assumed to be commutative and possess a unit.}

\section{Weight of a cohomology class with respect to a fibration}\label{secweight}

We  start by recalling a result from \cite{FG06} which will be used in this paper.

Let
$p\co E\to B$ be a fibration. The {\em Schwarz genus} of $p$,
denoted $\genus(p)$, is defined to be the least integer $k$ such that
the base $B$ may be covered by $k$ open sets $U_1,\ldots ,U_k$ on each
of which $p$ admits a local section (a continuous map $s_i\co
U_i\to E$ satisfying $p\circ s_i=\mathbf{1}_{U_i}$ for $i=1, \dots, k$).
The concept of genus of a fibration was introduced and thoroughly
studied by A.\ S.\ Schwarz
\cite{Sch66}. In the literature the term sectional category is also used.

There are two important special cases when one is motivated to study the Schwarz genus.
Let $X$ be a topological space. Consider the Serre path fibration $\pi_0\co P_0X\to X$ where
the total space $P_0X$ consists of all paths $\gamma$ in $X$
with $\gamma(0)$ equal to a fixed base point $x_0\in X$.
The projection $\pi_0$ takes a path $\gamma$ to $\gamma(1)\in X$. The genus of this fibration equals the Lusternik - Schnirelmann category $$\cat(X)=\genus(\pi_0).$$

Another special case is used to estimate complexity of the robot motion planning problem.
Let $X$ be a topological space and let $X^I$ be the space of all paths in $X$ equipped with the compact-open topology. The map $\pi\co X^I\to X\times X$ associates to a path $\gamma\in X^I$ the pair of its end points
$(\gamma(0),\gamma(1))$. It is a fibration and its genus is called the topological complexity of $X$, denoted $$\TC(X)=\genus(\pi).$$ One thinks of $X$ as being the configuration space of a mechanical system; then $\TC(X)$ measures the \lq\lq navigational complexity\rq\rq\, of $X$; see \cite{Far06} for more detail.

A well-known lower bound for the genus of a fibration $p\co E\to B$
is given by the cup-length of the kernel of the induced map $p^*$ in
cohomology. Taking coefficients in an arbitrary ring $R$, suppose there are classes $u_1,\ldots ,u_\ell\in H^*(B;R)$ with $$p^*(u_1)=\ldots =p^*(u_\ell)=0\in H^\ast(E;R)$$ such that
their cup product $u_1\cdots u_\ell\not= 0 \in H^\ast(B;R)$ is non-zero. Then one has
$\genus(p)>\ell,$ see \cite{Sch66}, Theorem 4.

The notion of category weight of E.\ Fadell and S.\ Husseini \cite{FH92}
generalises to an arbitrary fibration $p\co E\to B$ as follows, see \cite{FG06}.
\begin{definition}\label{def1}
Let $u\in H^*(B;G)$ be a cohomology class, where $G$ is an abelian group.
The  weight of $u$ with respect to $p:E\to B$, denoted $\wgt_p(u)$, is defined to be the
largest integer $k$ such that $f^*(u)=0\in H^\ast(Y;G)$ for all maps $f\co Y\to B$
with $\genus(f^\ast p)\leq k$.
\end{definition}
Here $f^\ast p:E'\to Y$ denotes the pull-back fibration
of $p$ along $f$ and the inequality $\genus(f^\ast p)\le k$ means that there exists an open cover $U_1\cup \dots \cup U_k=Y$ and continuous maps $\phi_i:U_i\to E$ such that $p\circ \phi_i = f|U_i$ for $i=1, \dots, k$.

Clearly $\wgt_p(u)\geq 0$ for all classes $u$, and $\wgt_p(u)\geq 1$
if and only if $p^*(u)=0$.

It is convenient to define the weight of the zero cohomology class as being $+\infty$.
\begin{prop}\label{weight1}
Let $u_i\in H^{d_i}(B;G_i)$ be cohomology classes, $i=1, \dots, \ell$, such that their cup-product
$u_1\cdots u_\ell\in H^d(B;G)$ is non-zero, where $d=d_1+\dots+d_\ell$ and $G=G_1\otimes \dots\otimes G_\ell$. Then
\begin{eqnarray}\label{genusgreater}
\genus(p)>\sum_{i=1}^\ell \wgt_p(u_i).
\end{eqnarray}
\end{prop}
For a proof see \cite{FG06}. The lower bound for $\genus(p)$ given by this Proposition
may improve upon that given by the cup-length of $\mathrm{ker}\,p^*$, if
we can find indecomposables $u\in H^*(B)$ with $\wgt_p(u)\geq 2$. In the next section we show how one may find such cohomology classes using cohomology operations.

\section{$\TC$-weights of cohomology classes and cohomology operations}\label{sectcweight}

In this section we study weights of cohomology classes in the context of
topological complexity.

\begin{definition}\label{def2}
The {\em $\TC$-weight} of a cohomology class
$u\in H^*(X\times X;G)$ is defined as its weight $\wgt_\pi(u)$, in the sense of Definition \ref{def1}, with respect to
the path fibration $\pi: X^I\to X\times X$.
\end{definition}

As in general, $\wgt_\pi(u)\geq 1$ if and only if $\pi^\ast(u)=0$. The latter condition can be replaced by $\Delta^\ast(u)=0$ where $\Delta: X\to X\times X$ is the diagonal. If the group of coefficients $G$ is a field then $u$ can be viewed as an element of $H^\ast(X;G)\otimes H^\ast(X;G)$ and the property
$\Delta^\ast(u)=0$ can be expressed by saying that $u$ is a
{\it zero-divisor} in the sense of \cite{Far03}.

To find classes with $\wgt_\pi(u)\geq 2$ we will need the following Lemma.

\begin{lemma}\label{lm4}
Let $f=(\varphi,\psi)\co Y\to X\times X$ be a map where $\varphi$,
$\psi$ denote the projections of $f$ onto the first and second factors
of $X\times X$, respectively. Then the genus $\genus(f^\ast\pi)$ of the induced fibration of $\pi:X^I\to X\times X$ along $f$ is less than or equal to $2$
if and only if $Y=A\cup B$, where $A$ and $B$ are open in $Y$ and
$\varphi|_A\simeq\psi|_A: A\to X$, $\varphi|_B\simeq\psi|_B:B\to X$.
\end{lemma}
\begin{proof}
We have a pull-back diagram
\begin{diagram}[height=2em]
f^\ast (X^I) & \rTo & X^I \\
\dTo<{f^\ast \pi}& & \dTo>{\pi} \\
Y &\rTo^{f} & X\times X.
\end{diagram}
The conclusion of the lemma follows immediately from the the following
statement: {\em There exists a local section of $f^\ast\pi$ over an open
  subset $A\subseteq Y$ if and only if $\varphi|_A\simeq \psi|_A$.}
We first remark that a local section $s_A\co A\to f^\ast (X^I)$ is the
same as a map $\mathscr{S}_A\co A\to X^I$ satisfying
$\pi\circ\mathscr{S}_A=f|_A$. Assume that such a map $\mathscr{S}_A$
exists. Then we may define a homotopy $F\co A\times I\to X$ from
$\varphi|_A$ to $\psi|_A$ by $$F(a,t)=\mathscr{S}_A(a)(t),\quad a\in
A,\, t\in I.$$
Conversely, suppose we have a homotopy $G\co A\times I\to X$ from
$\varphi|_A$ to $\psi|_A$. Then we may define our map
$\mathscr{S}_A\co A\to X^I$ by the formula
$$\mathscr{S}_A(a)=(t\mapsto G(a,t)),\quad a\in A,\, t\in I.$$
\end{proof}

We now describe a method for finding indecomposable classes with
$\TC$-weight more than one, using
cohomology operations.

 Let $R$ and $S$ be abelian groups. A {\em stable cohomology operation of degree $i$}
 \begin{eqnarray}\label{theta}\theta: H^\ast(-;R) \to H^{\ast +i}(-;S)
 \end{eqnarray}
is a family of natural transformations
$\theta: H^{n}(-;R)\to H^{n+i}(-;S)$, one for each $n\in\Z$, which
commute with the suspension isomorphisms, see \cite{MT}.
It follows that $\theta$
commutes with all Mayer-Vietoris connecting homomorphisms, and each homomorphism (\ref{theta}) is additive, i.e. is a group homomorphism.

\begin{definition}\label{defexcess}
The excess of a stable cohomology operation $\theta$, denoted $e(\theta)$,
is defined to be the largest integer $n$ such that $\theta(u)= 0$ for all cohomology classes $u\in H^m(X;R)$ with
$m<n$.
\end{definition}
Consider a few examples. For any extension $0\to R'\to R\to R''\to 0$ of abelian groups the Bockstein homomorphism $$\beta: H^n(-;R'') \to H^{n+1}(-;R')$$
has excess one. The excess of the Steenrod square $$Sq^i: H^\ast(-;\Z_2) \to H^{\ast+i}(-; \Z_2)$$
equals $i$ and for any odd prime $p$ the excess of the Steenrod power operation
$$P^i: H^n(-;\Z_p) \to H^{n+2i(p-1)}(-;\Z_p)$$
equals $2i$, see \cite{Hat02}, pages 489 - 490. More generally, the excess of a composition of Steenrod squares
$\theta=Sq^I = Sq^{i_1}Sq^{i_{2}} \dots Sq^{i_n}$ satisfies
$$e(\theta)\geq \max_{1\le k\le n}\{i_k-i_{k+1}-\dots -i_n\}.$$
It is easy to see that for an {\it admissible} sequence $I= i_1i_{2}\dots i_n$ (i.e. such that $i_k\ge 2\cdot i_{k+1}$ for all $k$)
the excess equals
$$e(\theta)= \sum_k (i_k-2i_{k+1}),$$ which coincides with the standard notion of excess, see \cite{MT}, page 27.

Any cohomology class $u\in H^j(X;R)$ determines a class
$$\overline{u}=1\times u-u\times 1\in H^j(X\times X;R)$$ where
$\times$ denotes the cohomology cross product. Note that
$\overline{u}$ is a zero-divisor and hence $\wgt_\pi(\overline{u})\geq 1$.
 Observe that
 $$\theta(\overline{u})=\theta(p_2^*(u)-p_1^*(u))=p_2^*(\theta(u))-p_1^*(\theta(u))=\overline{\theta(u)},$$
 by the naturality and additivity of $\theta$ (here $p_1,p_2\co X\times X\to X$ are the projections onto each factor).

 Our main result in this paper is:

\begin{thm}\label{thm1} Let $\theta: H^\ast(-;R) \to H^{\ast+i}(-;S)$ be a stable cohomology operation of degree $i$ and excess $e(\theta)\geq n$.
 Then for any cohomology class $u\in H^n(X;R)$ of dimension $n$ the class
 $\theta(\overline u) = \overline{\theta(u)}= 1\times \theta(u) - \theta(u)\times 1\in H^{n+i}(X\times X;S)$ has $\TC$-weight at least $2$. In symbols, \begin{eqnarray}
 \wgt_{\pi}(\overline{\theta(u)}))\geq 2.\end{eqnarray}
\end{thm}
\begin{proof}
Let $f=(\varphi,\psi)\co Y\to X\times X$ be a map with
$\genus(f^\ast \pi)\leq 2$. Then by Lemma \ref{lm4} one has $Y=A\cup B$ with
restrictions $\varphi|A\simeq\psi|A$ and $\varphi|B\simeq\psi|B$ being homotopic. Consider the
element $$f^*(\overline{u})=\psi^*(u)-\varphi^*(u)\in H^{n}(Y;R).$$ By the
Mayer-Vietoris sequence for $Y$,
$$\ldots\to H^{n-1}(A\cap B;R)\stackrel{\delta}{\to}H^n(Y;R)\to
H^n(A;R)\oplus H^n(B;R)\to\ldots$$
we have $f^*(\overline{u})=\delta(w)$ for some $w\in H^{n-1}(A\cap B;R)$. Hence,
$$f^\ast(\overline{\theta(u))} = f^*(\theta(\overline{u}))=\theta(f^*(\overline{u}))=\theta(\delta(w))=\delta(\theta(w))=0,$$
since $\theta$ is a stable operation of excess $\geq n$ and $w$ has degree $n-1$.
\end{proof}
Note that similar results holds in a more general situation when $\theta\co E^*\to F^{*+i}$
is a stable cohomology operation between extraordinary
cohomology theories.

\section{Motion planning in fibre spaces}\label{secfibre}

In this section we give an upper bound for the topological complexity of fibre spaces in terms of invariants of the base and fibre. It will be used in the following section in the study of lens spaces.

\begin{lemma}\label{fibre}
Let $p: E\to B$ be a Hurewicz fibration with fibre $F$. Then
\begin{eqnarray}\label{one}
\TC(E) \, \leq\,   \TC(F) \cdot\cat(B\times B).
\end{eqnarray}
\end{lemma}
\begin{proof} Denote $k=\cat(B\times B)$ and $\ell=\TC(F)$. Suppose that
$$B\times B = U_1\cup \dots \cup U_k, \quad F\times F =V_1\cup \dots \cup V_\ell$$
are open covers such that each inclusion $U_i\to B\times B$ is null-homotopic and there exists a continuous section
$s_i: V_i\to F^I$ of the end-point map $F^I\to F\times F$ for each $i=1, \dots, \ell$.
Fix a homotopy $h_j: U_j\to (B\times B)^I=B^I\times B^I$ connecting the inclusion $U_j\to B\times B$
with the constant map onto $(x_0, x_0)$. For $(x,y)\in U_j$ the image $h_j(x,y)$ is a pair of paths $(\alpha_{x,y}, \beta_{x,y})$ in $B$
satisfying $\alpha_{x,y}(0)=x$, $\alpha_{x,y}(1)=x_0$ and  $\beta_{x,y}(0)=y$, $\beta_{x,y}(1)=x_0$.

As in Chapter 2, \S 7 of Spanier \cite{Sp},  denote $\bar B=\{(e, \omega)\in E\times B^I; \omega(0)=p(e)\}$ and consider a lifting function
$\lambda: \bar B\to E^I$ where for $(e, \omega)\in \bar B$ the image $\lambda(e, \omega)$ is a path in $E$ covering $\omega$ which starts at $e$.

For $(e, e')\in (p\times p)^{-1}(U_j)$ consider $x=p(e)\in B$ and $y=p(e')\in B$ and the paths $\lambda(e, \alpha_{x, y})\in E^I$ and
$\lambda(e', \beta_{x, y})\in E^I$. The end points of these paths $a=\lambda(e, \alpha_{x, y})(1)$ and $b=\lambda(e', \beta_{x, y})(1)$ lie in the fibre
$F$ above $x_0$. This defines a continuous map
\begin{eqnarray}
k_j: (p\times p)^{-1}(U_j) \to F\times F. \end{eqnarray}
Now, we denote by $W_{j, i}\subset E\times E$ the preimage $k_j^{-1}(V_i)$, where $j=1, \dots, k$ and $i=1, \dots, \ell$.

It is clear that the family $\{W_{j, i}\}$ is an open cover of $E\times E$ and over each set $W_{j,i}$ there is a continuous section of $E^I\to E\times E$:
if $(e, e')\in W_{j,i}$ then the connecting path is concatenation of $\lambda(e, \alpha_{x, y})$, path $s_i(a, b)$
in the fibre $F$ connecting $a$ to $b$, and the reverse path to $\lambda(e', \beta_{x, y})$. Hence, $\TC(E)\leq k\ell$.
\end{proof}

We mention the following special cases:

\begin{cor}
Let $E$ be the total space of a fibration with fibre $F$ such that the base $B$ is homotopy equivalent to a sphere $S^k$. Then $\TC(E)\leq 3\cdot \TC(F)$.
\end{cor}

\begin{cor}
Let $E$ be the total space of a fibration with base $B$ and fibre $S^k$ where $k$ is odd. Then one has $\TC(E) \le 2\cdot \cat(B\times B)$.
\end{cor}
\noindent
{\bf Question:} Can one replace $\cat(B\times B)$ in Lemma \ref{fibre} by the potentially smaller number $\TC(B)$?
In other words, we ask if the following inequality $$\TC(E) \leq \TC(B) \cdot \TC(F)$$ holds for any fibration $F\to E\to B$. If the topological complexity $\TC(K^2)$ of the Klein bottle $K^2$ equals 5 (we do not know if it is indeed the case) it would provide a counterexample since $K^2$ fibers over $S^1$ with fibre $S^1$ and $\TC(S^1)=2$.

\section{Topological complexity of lens spaces}\label{seclens}

In this section we apply Theorem \ref{thm1} to the problem of computing topological
complexity of lens spaces. Let $m\ge 2$ be an integer. We regard the cyclic group $\Z_m$ as the
multiplicative group $\{1,\omega,\ldots ,\omega^{m-1}\}\subseteq\mathbf{C}$ of $m$-th roots of unity. This acts freely on the unit sphere $S^{2n+1}\subseteq \mathbf{C}^{n+1}$ by pointwise multiplication. The quotient is the lens space $$L_m^{2n+1}=S^{2n+1}/\Z_m.$$ In the literature this space is known as $L_m(1, 1, \dots, 1)$, see page 144 of \cite{Hat02}.

We start by improving a general upper bound $\TC(X)\leq 2\cdot \dim(X)+1,$ see \cite{Far03}, for the lens spaces. There are many known examples when
$\TC(X)=2\dim(X)+1$; however, for real projective spaces one has a better upper bound $\TC(\mathbf {RP}^n)\leq 2n$ which is an equality if and only if $n$ is a power of two, see \cite{FTY03}, \cite{Far06}.

\begin{cor}\label{corlens} For the topological complexity of lens spaces $L_m^{2n+1}$ one has
\begin{eqnarray}\TC(L_m^{2n+1}) \le 2\cdot \dim (L_m^{2n+1})=4n+2.
\end{eqnarray}
\end{cor}
\begin{proof}
There is a locally trivial fibration $S^1\to L_m^{2n+1}\to {\mathbf {CP}}^n$. Indeed, a point of $L_m^{2n+1}$ is an orbit of a cyclic group acting linearly on
 $S^{2n+1}$ and associating to any such orbit the corresponding complex line gives the fibration mentioned above.
 It is easy to see that $\TC(\mathbf {CP}^n)= \cat(\mathbf {CP}^n\times \mathbf {CP}^n)=2n+1$. Applying Lemma \ref{fibre} we find
\begin{eqnarray}
\TC(L_m^{2n+1}) \leq (2n+1) \times 2= 4n+2.
\end{eqnarray}
\end{proof}

Next we describe lower bounds for $\TC(L_m^{2n+1})$ based on Theorem \ref{thm1}.

\begin{thm}\label{lower} The topological complexity of the lens space $L_m^{2n+1}$ satisfies
\begin{eqnarray}\label{tclens2}
\TC(L_m^{2n+1})\geq 2\cdot (k+\ell )+2
\end{eqnarray}
for any pair of integers $k, \ell$ such that
 $m$ does not divide $\binom{k+\ell}{k}$ and $0\leq k, \ell \leq n$.
\end{thm}
\begin{proof} The cohomology $H^i(L_m^{2n+1};\Z_m)$ is $\Z_m$ for $0\le i\le 2n+1$ and vanishes for $i>2n+1$, see \cite{Hat02}.
As generators one can choose $x\in H^1(L_m^{2n+1};\Z_m)$ and
\begin{eqnarray}y=\beta(x) \in H^2(L_m^{2n+1};\Z_m),
\end{eqnarray}
where $\beta\co H^1(-;\Z_m)\to
H^2(-;\Z_m)$ is the mod $m$ Bockstein homomorphism, and then $H^\ast(L_m^{2n+1};\Z_m)$, as a graded algebra, coincides with
the factor-ring $\Z_m[x, y]/I$. Here $I$ is the ideal generated by
$y^{n+1}$ and $x^2-ay$
where $a\in \Z$ is given by
$$a=\left\{
\begin{array}{l}
0,\quad \mbox{if $m$ is odd},\\

m/2, \quad \mbox{if $m$ is even}
\end{array}
\right.
$$
(see \cite{Hat02}, Example 3E.2).

The K\"unneth Theorem gives $$H^*(L_m^{2n+1}\times L_m^{2n+1};\Z_m)\cong H^*(L_m^{2n+1};\Z_m)\otimes
H^*(L_m^{2n+1};\Z_m)$$
(see \cite{Hat02}, Theorem 3.16 where we take $R=\Z_m$). Therefore classes of the form $x^{s_1}y^{r_1}\otimes x^{s_2}y^{r_2}$,
where $s_i\in \{0,1\}$ and $r_i\in \{0, \dots, n\}$, $i=1,2$, form an additive basis of $H^\ast(L_m^{2n+1}\times L_m^{2n+1};\Z_m)$ viewed as a free $\Z_m$-module.

Since $\beta$ is a stable cohomology operation of excess 1, we have by Theorem \ref{thm1}
\begin{eqnarray}\label{weight}
\wgt_\pi(\overline{\beta(x)}))=\wgt_\pi(\overline{y})\geq 2,
\end{eqnarray}
where $ \overline y = 1\otimes y-y\otimes 1\in H^2(L_m^{2n+1}\times L_m^{2n+1};\Z_m)$ is a zero-divisor. If for some $0\le k, \ell \le n$ the binomial coefficient $\binom{k+\ell}{k}$
is not divisible by $m$ then the power $(\overline y)^{k+\ell}$ is nonzero since it contains the term $(-1)^k \binom{k+\ell}{k} y^k\otimes y^\ell$. 
The product $\bar x\cdot(\bar y)^{k+l}$ is also nonzero (for obvious reasons) and hence applying Proposition \ref{weight1}, we obtain $\TC(L_m^{2n+1}) \geq 2(k+\ell)+2$. This completes the proof.
\end{proof}

To state the following result we need a new notation.
For an integer $n$ we will denote by $\alpha(n)=\alpha_2(n)$ the number of ones in the dyadic expansion of $n$.
To define a similar number $\alpha_p(n)$, for any odd prime $p$, consider the $p$-adic representation of $n$,
$$n=n_0+n_1p+\dots+n_kp^k, \quad n_i\in \{0,1,\dots, p-1\}.$$
The number $\alpha_p(n)$ is defined by counting indices $i$ such that $2n_i\geq p$, but our counting involves certain multiplicities $r_i=r_i(n)$. We set $r_i=0$ iff $2n_i<p$.
If $2n_i\geq p$, we denote by $r_i\geq 1$ the maximal $r\geq 1$ such that $n_{i+1}=n_{i+2}=\dots=n_{i+r-1}=(p-1)/2$. Thus $r_i=1$ iff $2n_i\geq p$ and $n_{i+1}\not= (p-1)/2$. Finally we define
\begin{eqnarray}
\alpha_p(n)=\sum_{i\geq 0} r_i.
\end{eqnarray}
Examples: $\alpha_3(13)=0$, $\alpha_3(14)=3$.

\begin{thm}\label{thmlens}  The topological complexity of the lens space $L_m^{2n+1}$ equals
\begin{eqnarray}\label{tclens}
\TC(L_m^{2n+1})=2\cdot \dim(L_m^{2n+1})=4n+2,
\end{eqnarray} assuming that $m$ is divisible by $p^{\alpha_p(n)+1}$,  for some prime $p$.
\end{thm}

\begin{proof} By Lemma \ref{div} in Appendix B, the maximal power of $p$ dividing $\binom{2n}{n}$ is
$p^{\alpha_p(n)}$. Hence, the assumption of Theorem \ref{thmlens} can be equivalently expressed by saying that
$m$ does not divide $\binom{2n}{n}$. The result follows by combining the upper bound of Corollary \ref{corlens}
with the lower bound given by Theorem \ref{lower} (where we take $k=\ell=n$). \end{proof}

The following statement is a useful special case of the previous theorem:

\begin{thm} \label{coradic} Suppose that $p$ is an odd prime and $n$ is such that its $p$-adic expansion,
$$n=n_0+n_1\cdot p+ \dots +n_k\cdot p^k,\quad \mbox{where}\quad n_i\in \{0, 1, \dots, m-1\},$$ involves only \lq\lq digits\rq\rq\, $n_i$
satisfying $n_i\le (p-1)/2.$ Then the topological complexity of the $(2n+1)$-dimensional lens space $L_p^{2n+1}$ equals
\begin{eqnarray}
\TC(L_p^{2n+1}) = 2\cdot \dim (L_p^{2n+1}) = 4n+2.\end{eqnarray}
\end{thm}
\begin{proof} The claim follows from Theorem \ref{thmlens} since the assumption of Theorem \ref{coradic} is equivalent to $\alpha_p(n)=0$.
\end{proof}

In the special case of $m=3$ Theorem \ref{coradic} applies and gives an explicit expression for topological complexity of lens spaces $L_3$ of dimensions $3, 7, 9, 19, 21, 25,  \dots$.

In the case $m=5$
Theorem \ref{coradic} applies to the lens spaces $L_5$ of dimensions $3,5,11,13,15, 21, \dots$.

\begin{thm} Assume that $m=2^r$. Then one has
\begin{eqnarray}
\TC(L_m^{2n+1}) = 2\cdot \dim(L_m^{2n+1})= 4n+2
\end{eqnarray}
for lens spaces $L_m$ of dimension $2n+1$ for all $n$ satisfying $\alpha(n) \le  r-1$ (i.e. for all $n$ which are sums of at most $r-1$ powers of $2$). \end{thm}
Recall that
$\alpha(n)$ denotes the number of ones in the dyadic expansion of $n$.

We see that the topological complexity of lens spaces $L_4^{2n+1}$ equals twice the dimension for all $n$ which are powers of $2$. The topological complexity of lens spaces $L_{8}^{2n+1}$ equals twice the dimension for all $n$ having at most two ones in their dyadic expansion; in other words, $n$ must be the sum of at most two powers of two.

\begin{cor}
The topological complexity of the $3$-dimensional lens space $L^3_m$ equals $6$ for all $m\geq 3$. On the other hand $\TC(L^3_2)=4$. 
\end{cor}
\begin{proof}
The above results imply that $\TC(L_m^3)=6$ for all $3$-dimensional lens spaces except possibly for $m=2$. In the remaining case one has
$\TC(L_2^3)=4$ as shown in \cite{FTY03} (note that $L_2^3={\bf {RP}}^3$).
\end{proof}

\begin{remark}{\em
Theorem \ref{thmlens} improves a result of J.\ Gonz\'alez
(\cite{Gon05}, Theorem 2.9) which states that $\TC(L_m^{2n+1})\geq 4n+1$ if
$m$ does not divide $\binom{2n}{n}$, and that if in addition $m$ is
even, then $\TC(L_m^{2n+1})$ equals either $4n+2$ or $4n+1$. Paper \cite{Gon06} of
J.\ Gonz\'alez
contains results concerning
$\TC(L_4^{2n+1})$ obtained using Brown-Peterson cohomology.
Papers \cite{Gon05}, \cite{Gon06} contain also a general discussion comparing the problem of computing the topological complexity of lens spaces and the immersion problem for lens spaces, inspired by the result of \cite{FTY03}.}
\end{remark}

\section*{Appendix A: Category weight of Fadell and Huseini}\label{secap}

In this appendix we give a short proof of a result in the spirit of theorems of  Fadell and
Husseini \cite{FH92}. Our Theorem \ref{FH} is slightly stronger than \cite{FH92} since it gives lower bounds for the strict category weight of Y. Rudyak \cite{Rud99} instead of the original category weight of \cite{FH92}.

\begin{definition}[Rudyak, \cite{Rud99}] {\rm
Let $u\in H^*(X)$ be a cohomology class. The {\em strict category weight} of $u$,
denoted $\swgt(u)$, is defined to be the largest integer $k$ such that
$f^*(u)=0$ for all maps $f\co Y\to X$ with $\cat(f)\leq k$. Recall that $\cat(f)\leq k$ means that
$Y$ may be
covered by open sets $U_1,\ldots ,U_k$, the restriction of $f$ to each
of which is null-homotopic.}
\end{definition}
Clearly $\swgt(u)$ coincides with the weight $\wgt_{\pi_0}(u)$ of $u$ with respect to the Serre fibration $\pi_0: P_0(X) \to X$.
One may improve on the classical cup-length lower bound by
finding indecomposable cohomology classes of category weight more than
one. The following result includes Theorem 3.12 of
\cite{FH92} as a special case, see also Corollary 4.7 of \cite{Rud99}.
\begin{thm}\label{FH}
Let $\theta$ be a stable cohomology operation $\theta: H^\ast(-;R)\to H^{\ast+i}(-;S)$ having excess\footnote{The notion of excess is described in Definition \ref{defexcess}.} $e(\theta)\geq n>0$, and let $u\in H^n(X;R)$ be a cohomology class of dimension
$n$. Then the class $\theta(u)\in H^{n+i}(X;S)$ has strict category weight greater than or equal to two,
\begin{eqnarray}\swgt(\theta(u))\geq 2.\end{eqnarray}
\end{thm}
\begin{proof}
We must show that $f^*(\theta(u))=0$ for all maps $f\co Y\to X$ with
$\cat(f)\leq 2$. Let $f$ be such a map. Then $Y=A\cup B$ where $A$ and
$B$ are open sets in $Y$ such that the restrictions $f|A$ and $f|B$ are
null-homotopic. By the Mayer-Vietoris sequence,
$$\ldots\to H^{n-1}(A\cap B)\stackrel{\delta}{\to}H^n(Y)\to
H^n(A)\oplus H^n(B)\to\ldots$$
we have that $f^*(u)=\delta (w)$ for some $w\in H^{n-1}(A\cap
B)$. Hence
$$f^*(\theta(u))=\theta(f^*(u))=\theta(\delta(w))=\delta(\theta(w))=0,$$
since $\theta$ has excess $\ge n$ and $w$ has degree $n-1$.
\end{proof}

\section*{Appendix B: Divisibility of binomial coefficients}

For convenience of the reader we include the following well-known result:

\begin{lemma}\label{div} Let $p$ be a prime and
$n=n_0+n_1p+n_2p^2 +\dots$ and $m=m_0+m_1p+m_2p^2 +\dots$
be $p$-adic representations of integers $n$ and $m$, where
$0\, \le\,  n_i, \, m_i\, <p$. The maximal integer $\ell$ such that $p^\ell$ divides the binomial coefficient $\binom{n+m}{n}$
equals the number of indices $i=0, 1, 2, \dots$ such that either
\begin{eqnarray}\label{geq}
n_i+m_i \geq p
\end{eqnarray}
or, for some $r\geq 1$, one has
\begin{eqnarray}\label{eq}
\begin{array}{l} n_i+m_i=n_{i-1}+m_{i-1}=\dots= n_{i-r}+m_{i-r}= p-1, \\ \\
n_{i-r-1}+m_{i-r-1}\geq p. \end{array}
\end{eqnarray}
\end{lemma}
\begin{proof}
One observes that $n!=p^\ell\cdot n'$ where $n'$ is an integer mutually prime to $p$ and $\ell= \sum_{i\geq 1}[n/p^i]$.
Hence the maximal power of $p$ dividing $\binom{n+m}{n}=\frac{(n+m)!}{n!m!}$ equals
\begin{eqnarray}\label{sum}
\qquad \quad \sum_{i\geq 1} \left( \left[\frac{n+m}{p^i}\right] - \left[\frac{n}{p^i}\right] - \left[\frac{m}{p^i}\right]\right) = \sum_{i\geq 1}\left( \left\{\frac{n}{p^i}\right\}+ \left\{\frac{m}{p^i}\right\}- \left\{\frac{n+m}{p^i}\right\}\right).
\end{eqnarray}
The symbols $[x]$ and $\{x\}$ denote integral and fractional parts of $x=[x]+\{x\}$ respectively.
In the sums (\ref{sum}) each term is zero or one and hence the value of the sum equals the number of ones. This implies Lemma \ref{div} since a term of (\ref{sum}) equals one if and only if either (\ref{geq}), or (\ref{eq})
hold for index $i-1$.
\end{proof}

\end{document}